\documentclass[12pt,a4paper]{amsart}
\usepackage{amssymb,amsfonts,amsthm,amsmath}
\usepackage{graphicx} 

\theoremstyle{plain}
\newtheorem{theorem}{Theorem}

\newtheorem{definition}{Definition}

\newtheorem{proposition}{Proposition}

\numberwithin{equation}{section} \numberwithin{theorem}{section}
\numberwithin{lemma}{section} \numberwithin{definition}{section}
\numberwithin{corollary}{section}
\numberwithin{proposition}{section} \textheight =24cm
\begin{document}
\title[Analogues of $q$-series]{Dirichlet series analogues of $q$-shifted factorial and the $q$-Kummer sum.}
\author{Geoffrey B Campbell}
\address{Mathematical Sciences Institute,
         The Australian National University,
         Canberra, ACT, 0200, Australia}

\email{Geoffrey.Campbell@anu.edu.au} 

\begin{abstract}
Recently, the concept of a $D$-analogue was introduced by the
author. This is a Dirichlet series analogue for the already known
and well researched hypergeometric $q$-series. we consider the
$D$-analogues of the $q$-binomial coefficients, and a $D$-analogue
of the $q$-Kummer (Bailey-Daum) sum.
\end{abstract}

\maketitle

\section{Introduction} \label{S:intro}

In a recent paper~\cite{gC2003a} the idea of a $D$-analogue was
introduced.  This is a Dirichlet series analogue for the already
known and well researched hypergeometric $q$-series, often called
the basic hypergeometric series.  The $q$-series is itself an
analogue for the ordinary hypergeometric series developed by Gauss
in the early nineteenth century. We note that the ordinary
hypergeometric series was introduced by Gauss~\cite{cG1813} in
1813, while the $q$-series analogue's are originally due to
Heine~\cite{eH1847,eH1878} in the mid 19th century.

Both types of hypergeometric series have been the topics of
far-reaching development and application throughout the 20th
century. The full range of applications would be exhaustive, but
one only needs cite such works as those of Andrews~(\cite{gA1975}
to \cite{gA1977}), Askey~(\cite{rA1975} to \cite{rA1987}) and
Baxter~(\cite{rB1980} to \cite{rB1982}) to affirm that the impact
of $q$-series has been great.

The author proposes that the newly introduced $D$-series analogue
should likewise turn out to be the subject of future development
and application.

A superficial comparison of the simplest $_{2}(symbol)_{1}$ three
types of analogue hypergeometric series is

  \begin{equation} \label{E:1.1}
 _{2}F_{1}(a,b;c;z) = \sum_{k=0}^{\infty}{\frac{(a)_{k}(b)_{k}}{k!
 (c)_{k}}z^{k}},
  \end{equation}

  \begin{equation} \label{E:1.2}
_{2}\phi_{1}(a,b;c;q,z)=\sum_{k=0}^{\infty}{\frac{(a;q)_{k}(b;q)_{k}}{(q;q)_{k}(c;q)_{k}}z^{k}},
  \end{equation}

  \begin{equation} \label{E:1.3}
 _{2}\Theta_{1}(a,b;c;\gamma,z)
=\sum_{k=1}^{\infty}\frac{\sigma_{-\gamma}(a;k)
\sigma_{-\gamma}(b;k)}{\sigma_{-\gamma}(c;k)}\frac{1}{k^{z}}.
  \end{equation}

\begin{flushleft}
In these, each kind of series has its own factorial analogue
function as follows,\end{flushleft}

\begin{equation} \label{E:1.4}
    (a)_{k}=
    \begin{cases}
      1,                       &\text{$k=0$;}\\
      (a)(a+1)...(a+k-1),      &\text{$k=1,2,...$;}
    \end{cases}
\end{equation}

\begin{equation} \label{E:1.5}
    (a;q)_{k}=
    \begin{cases}
      1,                               &\text{$k=0$;}\\
      (1-a)(1-aq)...(1-aq^{k-1}),      &\text{$k=1,2,...$;}
    \end{cases}
\end{equation}

\begin{equation} \label{E:1.6}
    \sigma_{-\gamma}(a;k)=
    \begin{cases}
      1,                               &\text{$k=1$;}\\
            \dfrac{\sigma_{-\gamma}(k)\sigma_{-\gamma}(k\prod_{p|k}{p})...\sigma_{-\gamma}(k\prod_{p|k}{p^{a-2}})}
      {\sigma_{-\gamma}(1)\sigma_{-\gamma}(\prod_{p|k}{p})...\sigma_{-\gamma}(\prod_{p|k}{p^{a-2}})},      &\text{$k=2,3,...$.}
    \end{cases}
\end{equation}

In the case of (\ref{E:1.4}) the similarity to a factorial, $k!$,
is clear, and $(a)_{k}$ is known as the shifted factorial. In the
case of (\ref{E:1.5}), $(a;q)_{k}$ is called the $q$-shifted
factorial function. In (\ref{E:1.6}) $a$ must be an integer; a
restriction that does not apply to (\ref{E:1.4}) nor
(\ref{E:1.5}). For (\ref{E:1.6}) the implied ``{$D$-shifted
factorial}'' in relation to the factorial is more disguised. If
$k$ has a prime decomposition $p_{1}^{a_{1}}p_{2}^{a_{2}}\cdots
p_{m}^{a_{m}}$, then some insight is gained from writing the
fraction of (\ref{E:1.6}) as

\begin{equation} \label{E:1.7}
\dfrac{\sigma_{-\gamma}(p_{1}^{a_{1}}p_{2}^{a_{2}}\cdots
p_{m}^{a_{m}})
       \sigma_{-\gamma}(p_{1}^{a_{1}+1}p_{2}^{a_{2}+1}\cdots p_{m}^{a_{m}+1})...
       \sigma_{-\gamma}(p_{1}^{a_{1}+a-1}p_{2}^{a_{2}+a-1}\cdots p_{m}^{a_{m}+a-1})}
{\sigma_{-\gamma}(p_{1}^{0}p_{2}^{0}\cdots p_{m}^{0})
       \sigma_{-\gamma}(p_{1}^{1}p_{2}^{1}\cdots p_{m}^{1})...
       \sigma_{-\gamma}(p_{1}^{a-1}p_{2}^{a-1}\cdots
       p_{m}^{a-1})},
\end{equation}

\begin{flushleft}
in which the function $\sigma_{n}(k)$ is the sum of $n$th powers
of the divisors of $k$, whence (see~\cite{mA1972}, \cite{tA1976},
\cite{gH1971} or \cite{rS1989})\end{flushleft}

 \begin{equation} \label{E:1.8}
 \sigma_{-\gamma}(k)= \sum_{d|k}{d^{-\gamma}} =k^{-\gamma} \sigma_{\gamma}(k)
 =\prod_{i=1}^{m}{\frac{(1-p_{i}^{-(a_{i}+1)\gamma})}{(1-p_{i}^{-\gamma})}}.
 \end{equation}

We see from this that (\ref{E:1.7}) can be rewritten as

\begin{equation} \label{E:1.9}
\dfrac{\sigma_{\gamma}(p_{1}^{a_{1}}p_{2}^{a_{2}}\cdots
p_{m}^{a_{m}})
       \sigma_{\gamma}(p_{1}^{a_{1}+1}p_{2}^{a_{2}+1}\cdots p_{m}^{a_{m}+1})...
       \sigma_{\gamma}(p_{1}^{a_{1}+a-1}p_{2}^{a_{2}+a-1}\cdots p_{m}^{a_{m}+a-1})}
{\sigma_{\gamma}(p_{1}^{0}p_{2}^{0}\cdots p_{m}^{0})
       \sigma_{\gamma}(p_{1}^{1}p_{2}^{1}\cdots p_{m}^{1})...
       \sigma_{\gamma}(p_{1}^{a-1}p_{2}^{a-1}\cdots
       p_{m}^{a-1})k^{a\gamma}};
\end{equation}
and an alternative form of (\ref{E:1.6}) is

\begin{equation} \label{E:1.10}
    \sigma_{-\gamma}(a;k)=
    \begin{cases}
      1,                               &\text{$k=1$;}\\
            \dfrac{\sigma_{\gamma}(k)\sigma_{\gamma}(k\prod_{p|k}{p})...\sigma_{\gamma}(k\prod_{p|k}{p^{a-2}})}
      {\sigma_{\gamma}(1)\sigma_{\gamma}(\prod_{p|k}{p})...\sigma_{\gamma}(\prod_{p|k}{p^{a-2}})k^{a\gamma}},      &\text{$k=2,3,...$;}
    \end{cases}
\end{equation}
so then

 \begin{equation} \label{E:1.11}
 \sigma_{-\gamma}(a;k)= \dfrac{\sigma_{\gamma}(a;k)}{k^{a\gamma}}.
 \end{equation}

A further comparison between the three kinds of analogue is shown
by the three Gauss hypergeometric sum formulae in terms of gamma
functions, $q$-shifted gamma functions as given in
Askey~\cite{rA1978,rA1980,rA1987} for example, and (perhaps)
$D$-shifted gamma functions.  These are, for various conditions
given in respectively, Bailey~\cite[pages~2--3]{wB1935}, Gasper
and Rahman~\cite[pages~9--11]{gG1990}, and the author's
paper~\cite{gC2003a}, as

  \begin{equation} \label{E:1.12}
  _{2}F_{1}(a,b;c;1)=\frac{\Gamma(c)\Gamma(c-a-b)}{\Gamma(c-a)\Gamma(c-b)},\end{equation}

  \begin{equation} \label{E:1.13}
  _{2}\phi_{1}(a,b;c;q,c/ab)=\frac{(c/a;q)_{\infty}(c/b;q)_{\infty}}{(c;q)_{\infty}(c/ab;q)_{\infty}},
  \end{equation}

  \begin{equation} \label{E:1.14}
 _{2}\Theta_{1}(a,b;c;\gamma,(c-a-b)\gamma)
=\frac{\zeta(c;\gamma)_{\infty} \zeta(c-a-b;\gamma)_{\infty}}
{\zeta(c-a;\gamma)_{\infty} \zeta(c-b;\gamma)_{\infty}}.
  \end{equation}

In the latter we use the Riemann zeta function $\zeta(a)$ to
further define for positive integers $n$

  \begin{equation} \label{E:1.15}
 \zeta(a;\gamma)_{n}=\prod_{k=0}^{n-1}{\zeta\left((a+k)\gamma)\right)},
  \end{equation}

\begin{flushleft}
with the extension to $\zeta(a;\gamma)_{\infty}$ as $n \rightarrow
\infty$. Our notation for $\sigma_{-\gamma}(a;k)$ is suggestive of
the relation to a multidimensional divisor function. We shall
continue to explore this interpretation for
$\sigma_{-\gamma}(a;k)$. (\ref{E:1.14}) was given
in~\cite{gC2003a} for the conditions: positive integers $a$, $b$
and $c$ and $\Re\gamma\geq0$, with $\Re c\gamma$,
$\Re(c-a-b)\gamma$, $\Re(c-a)\gamma$, $\Re(c-a-b)\gamma$, each
$>1$. We furthermore, showed by examples, that many simple cases
of (\ref{E:1.14}) are easily accessible and new.
\end{flushleft}

In the present paper, as we did in \cite{gC2003a}, we derive two
separate classes of $D$-analogues. One of these involves the above
$\zeta(a;\gamma)_{n}$ function, and the other involves the Jordan
totient function $J_{n}(k)$ extended in a similar way to our
extending $\zeta(a)$ into $\zeta(a;\gamma)_{n}$. We give further
examples of the new class of Dirichlet series. In \cite{gC2003a}
we featured as examples of our new transform, the $D$-binomial
theorem, and the $D$-analogue of the $q$-Gauss $_{2}F_{1}$
summation formula (\ref{E:1.14}).

We next give $D$-analogues of the classical $q$-Kummer
(Bailey-Daum) summation formula. (For an account of this
$q$-series formulae see Gasper and Rahman~\cite[pages 14,
236]{gG1990}).

The Dirichlet series transformations from~\cite{gC2003a} are used
to obtain our new results. We apply two kinds of Euler product
operator: a) over all primes, and b) over the primes dividing
positive integer $m$. For the latter, we need the following set
$S_m$ defined as the restricted product results in a summation
such as (\ref{E:1.17}).

Assume the positive integer prime decomposition
$m=\prod_{i=1}^{t}{p_{i}^{a_{i}}}$, and that we associate a set
$S_{m}=\{ x \in Z^{+} : x=\prod_{i=1}^{t}{p_{i}^{b_{i}}} \text{
for each } b_{i} \text{ a non-negative integer}\}$ with this
number. Then for positive integers $n$ we have new $q$-binomial
theorem analogues

\begin{theorem} (Campbell~\cite{gC2003a})\label{T:1.1}
For positive integers $n$, $\Re\beta>1$, $\Re\gamma>0$,
  \begin{equation} \label{E:1.16}
\sum_{k=1}^{\infty}\frac{\sigma_{-\gamma}(n;k)}{k^{\beta}}=\prod_{k=0}^{n-1}{\zeta(\beta+k\gamma)},
  \end{equation}
  \begin{equation} \label{E:1.17}
\sum_{k \in
S_{m}}\frac{\sigma_{-\gamma}(n;k)\lambda(k)}{k^{\beta}}
=\prod_{k=0}^{n-1}\frac{1}{\sigma_{-(\beta+k\gamma)}(\prod_{p|m}{p})},
  \end{equation}
  \end{theorem}

\begin{flushleft}
where $\lambda(k)$  is the Liouville function,
  \begin{equation} \label{E:1.18}
\lambda(k)=\prod_{i=1}^{t}{(-1)^{a_{i}}} \text{ for each }
k=\prod_{i=1}^{t}{p_{i}^{a_{i}}}.
  \end{equation}

See Apostol~\cite[page 37]{tA1976}, or any of \cite{gH1971,
rS1989, eT1951} for classical accounts describing this function.
The Liouville function frequently arises in $D$-analogues where a
$q$-series with a specific negative parameter is
transformed.\end{flushleft}

Before stating the $D$-analogue for the Kummer theorem, we again
examine briefly the function $\sigma_{-\gamma}(n;k)$. Theorem
\ref{T:1.1} gives us a good starting point. In~\cite{gC2003a} we
stated the cases of (\ref{E:1.16}), letting successively, $n=1$,
$n=2$, $n=3$, $n=4$, with $\beta$ mapped onto $\beta + \gamma$,

\begin{equation}\label{E:1.19}
\sum_{k=1}^{\infty}\frac{\sigma_{-\gamma}(1;k)}{k^{\beta +
\gamma}}=\sum_{k=1}^{\infty}{\frac{1}{k^{\beta +
\gamma}}}=\zeta(\beta + \gamma),
\end{equation}

\begin{equation}\label{E:1.20}
\sum_{k=1}^{\infty}\frac{\sigma_{-\gamma}(2;k)}{k^{\beta +
\gamma}}=\sum_{k=1}^{\infty}{\frac{\sigma_{-\gamma}(k)}{k^{\beta +
\gamma}}}=\zeta(\beta + \gamma)\zeta(\beta+2\gamma),
\end{equation}

\begin{gather}
\sum_{k=1}^{\infty}\frac{\sigma_{-\gamma}(3;k)}{k^{\beta +
\gamma}}=\sum_{k=1}^{\infty}{\frac{\sigma_{-\gamma}(k)\sigma_{-\gamma}(k
\prod_{p|k}p)}
                   {\sigma_{-\gamma}(\prod_{p|k}p)} \frac{1}{k^{\beta + \gamma}}}\notag \\
                   =\zeta(\beta + \gamma)\zeta(\beta+2\gamma)\zeta(\beta+3\gamma), \label{E:1.21}
\end{gather}

\begin{gather}
\sum_{k=1}^{\infty}\frac{\sigma_{-\gamma}(4;k)}{k^{\beta +
\gamma}}=\sum_{k=1}^{\infty}{\frac{\sigma_{-\gamma}(k)\sigma_{-\gamma}(k\prod_{p|k}p)\sigma_{-\gamma}(k\prod_{p|k}p^{2})}
                   {\sigma_{-\gamma}(\prod_{p|k}p)\sigma_{-\gamma}(\prod_{p|k}p^{2})}
                   \frac{1}{k^{\beta + \gamma}}}\notag \\
                   =\zeta(\beta + \gamma)\zeta(\beta+2\gamma)\zeta(\beta+3\gamma)\zeta(\beta+4\gamma),
\label{E:1.22}
\end{gather}

and so on, valid where each of the Riemann zeta functions is in
its range of absolute convergence. It is clear from (\ref{E:1.19})
to (\ref{E:1.22}) that $\sigma_{-\gamma}(n;k)$ simplifies when $k$
is a prime, or a power of a prime. If $k$ is a product of distinct
primes, or an integer power of a product of distinct primes the
function is also simpler to examine. The following is quite easy
to ascertain from the above and from (\ref{E:1.6}).
\begin{theorem}\label{T:1.2}
For positive integers $a$ and for any prime $p$,
\begin{equation} \label{E:1.23}
    \sigma_{-\gamma}(a;p)=\sigma_{-\gamma}({p^{a-1}})= \dfrac{1-p^{-a\gamma}}{1-p^{-\gamma}}
\end{equation}
\end{theorem}

\begin{theorem}\label{T:1.3}
If $k$ is a product of distinct primes, then for positive integers
$a$,
\begin{equation} \label{E:1.24}
    \sigma_{-\gamma}(a;k)=\sigma_{-\gamma}({k^{a-1}})=
    \prod_{p|k}\dfrac{1-p^{-a\gamma}}{1-p^{-\gamma}}.
\end{equation}
\end{theorem}

\begin{theorem}\label{T:1.4}
If $k$ is the square of prime $p$, then for positive integers $a$,
\begin{equation} \label{E:1.25}
    \sigma_{-\gamma}(a;k)=
    \begin{cases}
      \sigma_{-\gamma}({p^{a}})=\dfrac{1-p^{-(a+1)\gamma}}{1-p^{-\gamma}},      &\text{$k=2$;}\\
      \\
      \dfrac{\sigma_{-\gamma}({p^{a}})}{\sigma_{-\gamma}({p})}
      =\dfrac{1-p^{-(a+1)\gamma}}{1-p^{-2\gamma}},      &\text{$k \neq 2.$}
    \end{cases}
\end{equation}
\end{theorem}

\begin{theorem}\label{T:1.5}
If $k$ is the square of a product of distinct primes, then for
positive integers $a$,
\begin{equation} \label{E:1.26}
    \sigma_{-\gamma}(a;k)=
    \begin{cases}
      \sigma_{-\gamma}({k^{a}})=\prod_{p|k}\dfrac{1-p^{-(a+1)\gamma}}{1-p^{-\gamma}},      &\text{$k=2$;}\\
      \\
      \dfrac{\sigma_{-\gamma}({k^{a}})}{\sigma_{-\gamma}({\sqrt{k}})}
      =\prod_{p|k}\dfrac{1-p^{-(a+1)\gamma}}{1-p^{-2\gamma}},      &\text{$k \neq 2.$}
    \end{cases}
\end{equation}
\end{theorem}

The general version of (\ref{E:1.19}) to (\ref{E:1.22}),

  \begin{equation} \label{E:1.27}
\sum_{k=1}^{\infty}\frac{\sigma_{-\gamma}(a;k)}{k^{\beta+\gamma}}=\prod_{k=1}^{a}{\zeta(\beta+k\gamma)},
  \end{equation}

  shows upon equating Dirichlet coefficients, that

\begin{theorem}\label{T:1.6}
If $a$ and $k$ are positive integers,
\begin{equation} \label{E:1.28}
    \sigma_{-\gamma}(a;k)=
          k^{-\gamma} \sum_{k=k_{1} k_{2} \cdots k_{a}}{\dfrac{1}{k_{1}^{\gamma} k_{2}^{2\gamma} \cdots
          k_{a}^{a\gamma}}}.
\end{equation}
\end{theorem}


\begin{center}
  \includegraphics[scale=0.75]{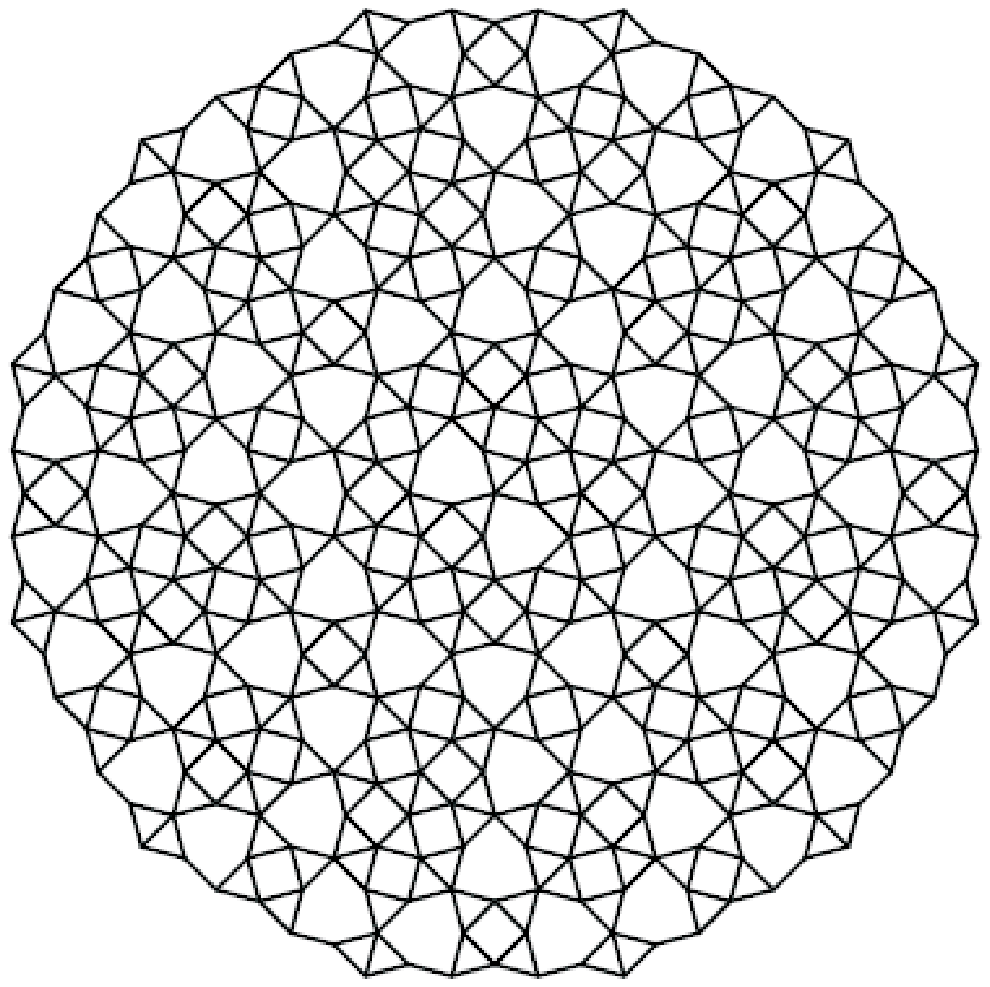}
\end{center}

\begin{center}
  figure 1.
\end{center}

This tiling pattern is enumerated recursively by a divisor
function related explicitly to the coefficient terms that occur in
the $D$-analogue of the binomial coefficients.

For a proof of this see Baake et al~\cite{mB2002a}. This
fascinating new connection between the theory of tilings arising
in the study of quasicrystals may be well worth further
investigation.

\section{Some further required notation.}\label{S:2}

In this section we restate and use the notation we defined for
$D$-analogues in~\cite{gC2003a}. That choice of nomenclature
resembled, where possible, the $q$-series notation. However it is
also our intention to display the similarity of our new
$D$-analogues with the ordinary hypergeometric series summations.
Firstly though, we use the notation given in Gasper and
Rahman~\cite{gG1990} for a product of $q$-shifted factorials:
\begin{definition}\label{D:2.1}
\begin{equation}\label{E:2.1}
(a_{1},a_{2},\ldots,a_{r};q)_{n}=(a_{1};q)_{n}(a_{2};q)_{n}\cdots(a_{r};q)_{n}.
\end{equation}
\end{definition}
Next we design the following definitions and notations. Firstly we
enlarge the definition given previously in (\ref{E:1.15}),

\begin{definition}\label{D:2.2}If $\zeta(a)$ is the Riemann zeta function and $\Re \gamma$
chosen such that the functions all exist, define for positive
integers $n$,
\begin{equation}\label{E:2.2}
\zeta(a;\gamma)_{n}=\prod_{k=0}^{n-1}{\zeta\left((a+k)\gamma)\right)}=\prod_{p}\frac{1}{(p^{-a\gamma};p^{-\gamma})_{n}},
\end{equation}
\begin{equation}\label{E:2.3}
\zeta(a_{1},a_{2},\ldots,a_{r};\gamma)_{n}=\zeta(a_{1};\gamma)_{n}\zeta(a_{2};\gamma)_{n}\cdots\zeta(a_{r};\gamma)_{n}.
\end{equation}\end{definition}

This definition is intended to bring a required notation for the
$D$-analogue of the normal gamma function $\Gamma(z)$ or its
$q$-analogue $\Gamma_q(z)$ found in the $q$-series literature such
as Gasper and Rahman~\cite{gG1990}. However, we shall have to
postpone our investigation of this feature to another paper.

\begin{definition}\label{D:2.3}If $\sigma_k(a)$ is the sum of $k$th powers of the divisors of positive integer~$a$  as in (\ref{E:1.6}) then for positive
integers $n$,
\begin{equation}\label{E:2.4}
\sigma_{-\gamma}(a;k)=\prod_{j=0}^{a-2}\frac{\sigma_{-\gamma}(k\prod_{p|k}p^{j})}{\sigma_{-\gamma}(\prod_{p|k}p^{j})},
\quad \text{defined as } 1 \text{ at } a=1, \text{ and}
\end{equation}
\begin{equation}\label{E:2.5}
\sigma_{-\gamma}(a_{1},a_{2},\ldots,a_{r};k)=\sigma_{-\gamma}(a_{1};k)\sigma_{-\gamma}(a_{2};k)\cdots\sigma_{-\gamma}(a_{r};k).
\end{equation}\end{definition}
With these definitions, we formed in~\cite{gC2003a} respectively,
the following theorems.

\begin{theorem} \label{T:2.1}
For positive integers $a_{i}$, $b_{i}$, $c_{i}$, $d_{i}$, with
$\Re z$ chosen for convergence,
 \begin{gather}
\prod _{p}{_{r+s}\phi_{r+s-1}\left[
\begin{matrix}p^{-a_{1}\gamma},p^{-a_{2}\gamma},\ldots,p^{-a_{r}\gamma},-p^{-c_{1}\gamma},-p^{-c_{2}\gamma},\ldots,-p^{-c_{s}\gamma}; & p^{-\gamma},p^{-z} \\
p^{-b_{1}\gamma},p^{-b_{2}\gamma},\ldots,p^{-b_{(r-1)}\gamma},-p^{-d_{1}\gamma},-p^{-d_{2}\gamma},\ldots,-p^{-d_{s}\gamma}
&
\end{matrix}
\right]}  \notag \\
=\sum_{k=1}^{\infty}{\frac{\sigma_{-\gamma}(a_{1},a_{2},\ldots,a_{r};k)\sigma_{-2\gamma}(c_{1},c_{2},\ldots,c_{s};k)
             \sigma_{-\gamma}(d_{1},d_{2},\ldots,d_{s};k)}
{\sigma_{-\gamma}(b_{1},b_{2},\ldots,b_{r-1};k)\sigma_{-2\gamma}(d_{1},d_{2},\ldots,d_{s};k)
 \sigma_{-\gamma}(c_{1},c_{2},\ldots,c_{s};k)}
\frac{1}{k^{z}}}, \quad and\label{E:2.6}
 \end{gather}

 \begin{gather}
\prod _{p|m}{_{r+s}\phi_{r+s-1}\left[
\begin{matrix}p^{-a_{1}\gamma},p^{-a_{2}\gamma},\ldots,p^{-a_{r}\gamma},-p^{-c_{1}\gamma},-p^{-c_{2}\gamma},\ldots,-p^{-c_{s}\gamma}; & p^{-\gamma},p^{-z} \\
p^{-b_{1}\gamma},p^{-b_{2}\gamma},\ldots,p^{-b_{(r-1)}\gamma},-p^{-d_{1}\gamma},-p^{-d_{2}\gamma},\ldots,-p^{-d_{s}\gamma}
&
\end{matrix}
\right]}  \notag \\
=\sum_{k \in
S_{m}}{\frac{\sigma_{-\gamma}(a_{1},a_{2},\ldots,a_{r};k)\sigma_{-2\gamma}(c_{1},c_{2},\ldots,c_{s};k)
             \sigma_{-\gamma}(d_{1},d_{2},\ldots,d_{s};k)}
{\sigma_{-\gamma}(b_{1},b_{2},\ldots,b_{r-1};k)\sigma_{-2\gamma}(d_{1},d_{2},\ldots,d_{s};k)
 \sigma_{-\gamma}(c_{1},c_{2},\ldots,c_{s};k)}
\frac{1}{k^{z}}}.\label{E:2.7}
 \end{gather}
\end{theorem}

Similarly, we have for $-p^{-z}$ in place of $p^{-z}$,

\begin{theorem} \label{T:2.2}
For positive integers $a_{i}$, $b_{i}$, $c_{i}$, $d_{i}$, with
$\Re z$ chosen for convergence,
 \begin{gather}
\prod _{p}{_{r+s}\phi_{r+s-1}\left[
\begin{matrix}p^{-a_{1}\gamma},p^{-a_{2}\gamma},\ldots,p^{-a_{r}\gamma},-p^{-c_{1}\gamma},-p^{-c_{2}\gamma},\ldots,-p^{-c_{s}\gamma}; & p^{-\gamma},-p^{-z} \\
p^{-b_{1}\gamma},p^{-b_{2}\gamma},\ldots,p^{-b_{(r-1)}\gamma},-p^{-d_{1}\gamma},-p^{-d_{2}\gamma},\ldots,-p^{-d_{s}\gamma}
&
\end{matrix}
\right]}  \notag \\
=\sum_{k=1}^{\infty}{\frac{\sigma_{-\gamma}(a_{1},a_{2},\ldots,a_{r};k)\sigma_{-2\gamma}(c_{1},c_{2},\ldots,c_{s};k)
             \sigma_{-\gamma}(d_{1},d_{2},\ldots,d_{s};k)}
{\sigma_{-\gamma}(b_{1},b_{2},\ldots,b_{r-1};k)\sigma_{-2\gamma}(d_{1},d_{2},\ldots,d_{s};k)
 \sigma_{-\gamma}(c_{1},c_{2},\ldots,c_{s};k)}
\frac{\lambda(k)}{k^{z}}}, \quad and\label{E:2.8}
 \end{gather}

 \begin{gather}
\prod _{p|m}{_{r+s}\phi_{r+s-1}\left[
\begin{matrix}p^{-a_{1}\gamma},p^{-a_{2}\gamma},\ldots,p^{-a_{r}\gamma},-p^{-c_{1}\gamma},-p^{-c_{2}\gamma},\ldots,-p^{-c_{s}\gamma}; & p^{-\gamma},-p^{-z} \\
p^{-b_{1}\gamma},p^{-b_{2}\gamma},\ldots,p^{-b_{(r-1)}\gamma},-p^{-d_{1}\gamma},-p^{-d_{2}\gamma},\ldots,-p^{-d_{s}\gamma}
&
\end{matrix}
\right]}  \notag \\
=\sum_{k \in
S_{m}}{\frac{\sigma_{-\gamma}(a_{1},a_{2},\ldots,a_{r};k)\sigma_{-2\gamma}(c_{1},c_{2},\ldots,c_{s};k)
             \sigma_{-\gamma}(d_{1},d_{2},\ldots,d_{s};k)}
{\sigma_{-\gamma}(b_{1},b_{2},\ldots,b_{r-1};k)\sigma_{-2\gamma}(d_{1},d_{2},\ldots,d_{s};k)
 \sigma_{-\gamma}(c_{1},c_{2},\ldots,c_{s};k)}
\frac{\lambda(k)}{k^{z}}}.\label{E:2.9}
 \end{gather}
\end{theorem}

Theorems \ref{T:2.1} and \ref{T:2.2} show how to take account of
negative parameters arising in the $q$-series under the product
operator in the context of this paper. Our notation always
requires equal numbers of  $c_{i}$  and $d_{i}$ terms in the
$q$-series, due to the necessary pairing of numerator and
denominator terms for cancellations. We next apply these theorems
to known basic hypergeometric series summations and transforms to
obtain new results as Dirichlet series analogues of the original
$q$-summation formulae. To do this we state the

\begin{definition}\label{D:2.8}
Either side of (\ref{E:2.6}), (\ref{E:2.7}), (\ref{E:2.8}) and
(\ref{E:2.9}) are defined respectively as:-
\begin{subequations}
 \begin{gather}
_{r+s}\Theta_{r+s-1}\left[
\begin{matrix}a_{1},a_{2},\ldots,a_{r},-\backslash c_{1},-\backslash c_{2},\ldots,-\backslash c_{s};    & \gamma,z \\
              b_{1},b_{2},\ldots,b_{r-1},-\backslash d_{1},-\backslash d_{2},\ldots,-\backslash d_{s} &
\end{matrix}
\right],\label{E:2.10a}  \\
_{r+s}\Theta_{r+s-1}\left[
\begin{matrix}m| & a_{1},a_{2},\ldots,a_{r},-\backslash c_{1},-\backslash c_{2},\ldots,-\backslash c_{s};    & \gamma,z \\
                 & b_{1},b_{2},\ldots,b_{r-1},-\backslash d_{1},-\backslash d_{2},\ldots,-\backslash d_{s} &
\end{matrix}
\right],\label{E:2.10b}  \\
_{r+s}\Theta_{r+s-1}\left[
\begin{matrix}a_{1},a_{2},\ldots,a_{r},-\backslash c_{1},-\backslash c_{2},\ldots,-\backslash c_{s};    & \gamma,-\backslash z \\
              b_{1},b_{2},\ldots,b_{r-1},-\backslash d_{1},-\backslash d_{2},\ldots,-\backslash d_{s} &
\end{matrix}
\right],\label{E:2.10c}  \\
_{r+s}\Theta_{r+s-1}\left[
\begin{matrix}m| & a_{1},a_{2},\ldots,a_{r},-\backslash c_{1},-\backslash c_{2},\ldots,-\backslash c_{s};    & \gamma,-\backslash z \\
                 & b_{1},b_{2},\ldots,b_{r-1},-\backslash d_{1},-\backslash d_{2},\ldots,-\backslash d_{s} &
\end{matrix}
\right],\label{E:2.10d}
 \end{gather}
 \end{subequations}
 where the $-\backslash$ preceding any variable denotes that it
 comes from a negative valued parameter in the $q$-series.
\end{definition}

In future use of the transforms of this section we will always
delineate between the negative and positive values of a parameter
even when this was previously not an issue of consideration in the
$q$-series itself.

\section{The transform applied to known $q$-series identities.}\label{S:3}

This section restates two theorems from Campbell~\cite{gC2003a} as
the basis of our results on the transformed $q$-series to
$D$-series. They show that a valid $q$-series identity or
transform maps onto a valid Dirichlet series identity or transform
in many cases. We shall frequently employ definition \ref{D:2.2}
here, and a counterpart to this using the Jordan totient function,
namely

\begin{definition}\label{D:3.1}Consider the Jordan totient function
 $J_{\gamma}(m)= m^{\gamma}\prod_{i=1}^{t}{(1-p_{i}^{-\gamma})},$ where $\gamma\geq0$, and  $m=\prod_{i=1}^{t}{p_{i}^{b_{i}}}$ is
the unique prime decomposition of $m$. Define then for positive
integers $m, n, a_{1}, a_{2}, \ldots, a_{r},$

\begin{equation}\label{E:3.1}
J(m| \ a_{1}, a_{2}, \ldots, a_{r};\gamma)_{n}
=\prod_{i=1}^{r}{\prod_{j=0}^{n-1}{\frac{J_{(a_{i}+j)\gamma}(m)}{m^{(a_{i}+j)\gamma}}}}
=\prod_{i=1}^{r}{\prod_{p|m}{(p^{-a_{i}\gamma};
p^{-\gamma})_{n}}}.
\end{equation}
\end{definition}

It is relatively easy to see from this and (\ref{E:1.8}) that
\begin{gather}
\frac{J(m| \ a_{1}, a_{2}, \ldots, a_{r};\gamma)_{n}}{J(m| \
b_{1}, b_{2}, \ldots, b_{r};\gamma)_{n}}
=\prod_{i=1}^{r}{\prod_{j=0}^{n-1}{\prod_{p|m}{\frac{(1-p^{-(a_{i}+j)\gamma})}{(1-p^{-(b_{i}+j)\gamma})}}}}
=\prod_{i=1}^{r}{\prod_{p|m}{\frac{(p^{-a_{i}\gamma};
p^{-\gamma})_{n}}{(p^{-b_{i}\gamma};
p^{-\gamma})_{n}}}}\notag \\
=\prod_{i=1}^{r}{\frac{\sigma_{-\gamma}(\prod_{p|m}{p^{a_{i}-1}})
\sigma_{-\gamma}(\prod_{p|m}{p^{a_{i}}})
\sigma_{-\gamma}(\prod_{p|m}{p^{a_{i}+1}}) \cdots
\sigma_{-\gamma}(\prod_{p|m}{p^{a_{1}+n-2}})}{
\sigma_{-\gamma}(\prod_{p|m}{p^{b_{i}-1}})
\sigma_{-\gamma}(\prod_{p|m}{p^{b_{i}}})
\sigma_{-\gamma}(\prod_{p|m}{p^{b_{i}+1}}) \cdots
\sigma_{-\gamma}(\prod_{p|m}{p^{b_{1}+n-2}})} },\label{E:3.2}
\end{gather}
and this will be used in simplifying the results from finite Euler
product transforms.

\begin{theorem}(see Campbell~\cite{gC2003a})\label{T:3.1}
Suppose that for each prime $p$ and positive integers $a_{i}$,
$b_{i}$, $c_{i}$, $d_{i}$, $e_{i}$, $f_{i}$, $g_{i}$, $h_{i}$, we
have a $q$-series identity of the generic form

\begin{subequations}
\begin{gather}
_{r+s}\phi_{r+s-1}\left[
\begin{matrix}p^{-a_{1}\gamma},p^{-a_{2}\gamma},\ldots,p^{-a_{r}\gamma},-p^{-c_{1}\gamma},-p^{-c_{2}\gamma},\ldots,-p^{-c_{s}\gamma}; & p^{-\gamma},p^{-z} \\
p^{-b_{1}\gamma},p^{-b_{2}\gamma},\ldots,p^{-b_{(r-1)}\gamma},-p^{-d_{1}\gamma},-p^{-d_{2}\gamma},\ldots,-p^{-d_{s}\gamma}
&
\end{matrix}
\right] \label{E:3.3a}\\
=\frac{(p^{-e_{1}\gamma},p^{-e_{2}\gamma},\ldots,p^{-e_{r}\gamma};p^{-\gamma})_{m_{1}}
\
(p^{-g_{1}\gamma},p^{-g_{2}\gamma},\ldots,p^{-g_{s}\gamma};p^{-2\gamma})_{m_{2}}\cdots}
{(p^{-f_{1}\gamma},p^{-f_{2}\gamma},\ldots,p^{-f_{r}\gamma};p^{-\gamma})_{n_{1}}
\
(p^{-h_{1}\gamma},p^{-h_{2}\gamma},\ldots,p^{-h_{s}\gamma};p^{-2\gamma})_{n_{2}}\cdots}.
\label{E:3.3b}
\end{gather}
\end{subequations}

Then the two Euler product transforms applied over primes to both
sides of this yield respectively, the two $D$-analogue summation
formulae

\begin{subequations}
\begin{gather}
_{r+s}\Theta_{r+s-1}\left[
\begin{matrix}a_{1},a_{2},\ldots,a_{r},-\backslash c_{1},-\backslash c_{2},\ldots,-\backslash c_{r};    & \gamma,z \\
              b_{1},b_{2},\ldots,b_{r-1},-\backslash d_{1},-\backslash d_{2},\ldots,-\backslash d_{r} &
\end{matrix}
\right] \label{E:3.4a}\\
=\frac{\zeta(f_{1},f_{2},\ldots,f_{s_{1}};\gamma)_{n_{1}} \
\zeta(h_{1},h_{2},\ldots,h_{s_{2}};2\gamma)_{n_{2}}\cdots}
      {\zeta(e_{1},e_{2},\ldots,e_{r_{1}};\gamma)_{m_{1}} \
      \zeta(g_{1},g_{2},\ldots,g_{r_{2}};2\gamma)_{m_{2}}\cdots}, \quad and  \label{E:3.4b}
\end{gather}
\end{subequations}
\begin{subequations}
\begin{gather}
_{r+s}\Theta_{r+s-1}\left[
\begin{matrix}m| & a_{1},a_{2},\ldots,a_{r},-\backslash c_{1},-\backslash c_{2},\ldots,-\backslash c_{r};    & \gamma,z \\
                 & b_{1},b_{2},\ldots,b_{r-1},-\backslash d_{1},-\backslash d_{2},\ldots,-\backslash d_{r} &
\end{matrix}
\right] \label{E:3.5a}\\
=\frac{J(m| \ e_{1}, e_{2}, \ldots, e_{r_{1}};\gamma)_{m_{1}} \
J(m| \ g_{1}, g_{2}, \ldots, g_{r_{1}};2\gamma)_{m_{2}}\cdots}
{J(m| \ f_{1}, f_{2}, \ldots, f_{r_{2}};\gamma)_{n_{1}} \ J(m| \
h_{1}, h_{2}, \ldots, h_{r_{2}};2\gamma)_{n_{2}}\cdots}.
\label{E:3.5b}
\end{gather}
\end{subequations}
\end{theorem}
A simple restatement of this with appropriate negative parameter
is
\begin{theorem}(see Campbell~\cite{gC2003a})\label{T:3.2}
Theorem \ref{T:3.1} with $-p^{-z}$ replacing $p^{-z}$ in
(\ref{E:3.3a}), together with $-\backslash z$ replacing $z$ in
each of (\ref{E:3.4a}) and (\ref{E:3.5a}) is true.
\end{theorem}
These are the required theorems to enable us to write down at a
glance, many $D$-series analogues of $q$-series. They assert no
claims as to the convergence of the resultant series of functions,
but formally give the $D$-analogues without deeper considerations.
It may, for example, turn out that the question of convergence of
the Euler product is non-trivial, and that various deeper analysis
will be required to justify resulting formulae. Not at this stage
concerning ourselves with the finer questions, we next show some
simple examples, of which all appear to be new results.

\section{The $D$-Kummer theorem.}\label{S:4}

In this section we apply the transforms in theorems \ref{T:3.1}
and \ref{T:3.2} to the so-called "Kummer" theorems in ordinary and
basic hypergeometric series. In later papers we may go into more
detail as to the significance of these new results. But for now we
shall content ourselves with statement of the summation formulae.

The ordinary Kummer theorem and its $q$-analogue theorem are shown
in Gasper and Rahman~\cite[page~14]{gG1990}) to be respectively,

\begin{equation}\label{E:4.1}
_{2}F_{1}(a, b; 1+a-b ; -1)
=\frac{\Gamma(1+a-b)\Gamma(1+{\frac{1}{2}}a)}
 {\Gamma(1+a) \Gamma(1+ {\frac{1}{2}}a - b)},
\end{equation}

\begin{equation}\label{E:4.2}
_{2}\phi_{1}(a, b;aq/b;q,-q/b)
=\frac{(-q;q)_{\infty} (aq,
aq^{2}/b^{2};q^{2})_{\infty}}{(aq/b, -q/b; q)_{\infty}}.
\end{equation}

For a proof of the ordinary hypergeometric series result
(\ref{E:4.1}) see Bailey~\cite[pages 9--10]{wB1935}. A
straightforward application of theorem \ref{E:3.2} to
(\ref{E:4.2}) yields the $D$-analogue formulae

\begin{proposition}\label{P:4.1}(The ``would-be'' $D$-Kummer Theorem)
If each of $a$ and $b$ are positive integers and $\gamma>0$ such
that each of the Riemann zeta functions in (\ref{E:4.3}) are
defined from convergent Euler products,
\begin{equation}\label{E:4.3} _{2}\Theta_{1}(a,
b;1+a-b;\gamma,-\backslash (1-b)\gamma) =\frac{\zeta(1, 1+a-b ;
\gamma)_{\infty} \zeta(1-b; 2\gamma)_{\infty}}
{\zeta(1-b;\gamma)_{\infty} \zeta(1+a,
1+{\frac{1}{2}}a-b;2\gamma)_{\infty}};
\end{equation}
and
\begin{equation}\label{E:4.4}
_{2}\Theta_{1}(m| \; a, b;1+a-b;\gamma,-\backslash (1-b)\gamma)
=\frac{J(m| \; 1-b;\gamma)_{\infty} J(m| \; 1+a,
1+{\frac{1}{2}}a-b;2\gamma)_{\infty}}{J(m| \; 1, 1+a-b ;
\gamma)_{\infty} J(m| \; 1-b; 2\gamma)_{\infty}} .
\end{equation}
\end{proposition}

 According to theorem \ref{E:2.2} applied to the left side of (\ref{E:4.2}), the left sides of (\ref{E:4.3})
 and (\ref{E:4.4}) are respectively

 \begin{gather}
\sum_{k=1}^{\infty}{\frac{\sigma_{-\gamma}(a,b;k)}
{\sigma_{-\gamma}(1+a-b;k)}
\frac{\lambda(k)}{k^{(1-b)\gamma}}}, \label{E:4.5} \quad and \\
\sum_{k \in S_{m}}{\frac{\sigma_{-\gamma}(a,b;k)}
{\sigma_{-\gamma}(1+a-b;k)}
\frac{\lambda(k)}{k^{(1-b)\gamma}}}.\label{E:4.6}
 \end{gather}

A quick inspection of the oscillating series (\ref{E:4.5})
indicates a general problem with its convergence. It is not at all
obvious whether there is any case of the series that converges
either absolutely or conditionally under the proposed conditions
for the proposition. The term $k^{-(1-b)\gamma}$ for $k>1$ with
$b$ a fixed positive integer does not approach zero as $k$
increases. This means we will need to rely upon the coefficient
terms in (\ref{E:4.5}) and  (\ref{E:4.6}) to achieve convergence.
Moreover, we know what the average order of $\sigma_{-\gamma}(k)$
is from the
\begin{theorem}\label{T:4.1}(see Apostol~\cite[pages~60-61]{tA1976})
If $\gamma>0$ and $x>1$,
\begin{equation}\label{E:4.7}
\sum_{k \leq x}{\sigma_{\gamma}}(k) =
\begin{cases}\frac{1}{2}\zeta(2)x^{2} + O(x log x), & \gamma=1;\\
\\
\dfrac{\zeta(\gamma + 1)}{\gamma + 1} x^{\gamma + 1} +
O(x^{\beta}), & \gamma\neq1, \; \beta = max\{1,\gamma\};
\end{cases}
\end{equation}
\begin{equation}\label{E:4.8}
\sum_{k \leq x}{\sigma_{-\gamma}}(k) =
\begin{cases}\zeta(2)x + O(log x), & \gamma=1;\\
\\
\zeta(\gamma + 1)x + O(x^{\beta}), & \gamma\neq1, \; \beta =
max\{0,1-\gamma\}.
\end{cases}
\end{equation}
\end{theorem}
Applying theorem \ref{T:4.1}, to the definitions of
$\sigma_{-\gamma}(a;k)$ of (\ref{E:1.6}) and (\ref{E:1.10}) we can
expect to obtain estimates of the behavior of the coefficients in
a $D$-series, such as for those presented in proposition
\ref{P:4.1}. Also, examining the coefficients of  (\ref{E:4.5})
using (\ref{E:1.23}) tells us, for instance, that if $k$ is a
prime
\begin{equation}\label{E:4.9}
\frac{\sigma_{-\gamma}(a,b;k)} {\sigma_{-\gamma}(1+a-b;k)}
= \prod_{p|k} \dfrac{(1-p^{-a\gamma})(1-p^{-b\gamma})}{(1-p^{-\gamma})(1-p^{-(1+a-b)\gamma})}.
\end{equation}
So, it seems \emph{per se} that the series  (\ref{E:4.5}) and
(\ref{E:4.6}) will diverge for any given fixed positive integer
$b$, for any positive integer $a$ chosen. However, this is not
strictly the case, as the problem with divergence in (\ref{E:4.5})
arises from the divergence of the Euler product over \emph{all}
primes. Since (\ref{E:4.6}) comes from a finite Euler product, it
is nonetheless true in Proposition~\ref{P:4.1}, and therefore the
$D$-series analogue (\ref{E:4.6}) of the Kummer theorems could on
its own perhaps be rated as a theorem. Let us state this
therefore, simplifying using (\ref{E:3.2}), to obtain our
consolation prize,

\begin{theorem}\label{T:4.2}(The ``half-done'' $D$-Kummer Theorem)
If each of $\frac{1}{2}a$ and $b$ are positive integers and
$\gamma>0$,
\begin{gather}\label{E:4.10}
_{2}\Theta_{1}(m| \; a, b;1+a-b;\gamma,-\backslash (1-b)\gamma)
=\frac{J(m| \; 1-b;\gamma)_{\infty} J(m| \; 1+a,
1+{\frac{1}{2}}a-b;2\gamma)_{\infty}}{J(m| \; 1, 1+a-b ;
\gamma)_{\infty} J(m| \; 1-b; 2\gamma)_{\infty}} \notag \\
= \prod_{j=0}^{\infty} {\frac{\sigma_{-\gamma}(\prod_{p|m}p^{j-b})
                         \sigma_{-2\gamma}(\prod_{p|m}p^{j+a})
                         \sigma_{-2\gamma}(\prod_{p|m}p^{j+\frac{1}{2}a-b})}
                        {\sigma_{-\gamma}(\prod_{p|m}p^{j})
                         \sigma_{-\gamma}(\prod_{p|m}p^{j+a-b})
                         \sigma_{-2\gamma}(\prod_{p|m}p^{j+b})}}.
\end{gather}
\end{theorem}

Obviously, the convergence questions of this section will
re-emerge in our sequel papers. Nevertheless it is clear that the
theorems in section~\ref{S:3} will still yield many
non-problematic $D$-analogues, such as the $D$-Dixon theorem, the
paper for which is in preparation. In the case of the $D$-Dixon
theorem there emerge many examples that are new, as was the case
for the $D$-Gauss analogue formula in~\cite{gC2003a} given in
(\ref{E:1.14}).

Finally it may serve our purpose to present the first several cases of the theorem, namely the first few integer values of $a$ and $b$ substituted. We can use (\ref{E:4.5}) and (\ref{E:4.6}) whence the following coefficient terms are apt.

\begin{equation}\label{E:4.11}
\frac{\sigma_{-\gamma}(a,b;k)} {\sigma_{-\gamma}(1+a-b;k)}
= \prod_{p|k} \dfrac{(1-p^{-a\gamma})(1-p^{-b\gamma})}{(1-p^{-\gamma})(1-p^{-(1+a-b)\gamma})}.
\end{equation}
\begin{equation}\label{E:4.12}
\frac{\sigma_{-\gamma}(2,1;k)} {\sigma_{-\gamma}(2;k)}
= \prod_{p|k} \dfrac{(1-p^{-2\gamma})(1-p^{-\gamma})}{(1-p^{-\gamma})(1-p^{-2\gamma})}=1.
\end{equation}
\begin{equation}\label{E:4.13}
\frac{\sigma_{-\gamma}(4,2;k)} {\sigma_{-\gamma}(3;k)}
= \prod_{p|k} \dfrac{(1-p^{-4\gamma})(1-p^{-2\gamma})}{(1-p^{-\gamma})(1-p^{-3\gamma})}.
\end{equation}
\begin{equation}\label{E:4.14}
\frac{\sigma_{-\gamma}(4,3;k)} {\sigma_{-\gamma}(2;k)}
= \prod_{p|k} \dfrac{(1-p^{-4\gamma})(1-p^{-3\gamma})}{(1-p^{-\gamma})(1-p^{-2\gamma})}.
\end{equation}
\begin{equation}\label{E:4.15}
\frac{\sigma_{-\gamma}(6,1;k)} {\sigma_{-\gamma}(6;k)}
= \prod_{p|k} \dfrac{(1-p^{-6\gamma})(1-p^{-\gamma})}{(1-p^{-\gamma})(1-p^{-6\gamma})}=1.
\end{equation}
\begin{equation}\label{E:4.16}
\frac{\sigma_{-\gamma}(6,2;k)} {\sigma_{-\gamma}(5;k)}
= \prod_{p|k} \dfrac{(1-p^{-6\gamma})(1-p^{-2\gamma})}{(1-p^{-\gamma})(1-p^{-5\gamma})}.
\end{equation} 

These coefficients when substituted give us the following particular cases of the $D$-Kummer theorem \ref{T:4.2} above.

\begin{equation}\label{E:4.17}
\frac{\sigma_{-\gamma}(a,b;k)} {\sigma_{-\gamma}(1+a-b;k)}
= \prod_{p|k} \dfrac{(1-p^{-a\gamma})(1-p^{-b\gamma})}{(1-p^{-\gamma})(1-p^{-(1+a-b)\gamma})}.
\end{equation}
\begin{equation}\label{E:4.18}
\frac{\sigma_{-\gamma}(2,1;k)} {\sigma_{-\gamma}(2;k)}
= \prod_{p|k} \dfrac{(1-p^{-2\gamma})(1-p^{-\gamma})}{(1-p^{-\gamma})(1-p^{-2\gamma})}=1.
\end{equation}
\begin{equation}\label{E:4.19}
\frac{\sigma_{-\gamma}(4,2;k)} {\sigma_{-\gamma}(3;k)}
= \prod_{p|k} \dfrac{(1-p^{-4\gamma})(1-p^{-2\gamma})}{(1-p^{-\gamma})(1-p^{-3\gamma})}.
\end{equation}
\begin{equation}\label{E:4.20}
\frac{\sigma_{-\gamma}(4,3;k)} {\sigma_{-\gamma}(2;k)}
= \prod_{p|k} \dfrac{(1-p^{-4\gamma})(1-p^{-3\gamma})}{(1-p^{-\gamma})(1-p^{-2\gamma})}.
\end{equation}
\begin{equation}\label{E:4.21}
\frac{\sigma_{-\gamma}(6,1;k)} {\sigma_{-\gamma}(6;k)}
= \prod_{p|k} \dfrac{(1-p^{-6\gamma})(1-p^{-\gamma})}{(1-p^{-\gamma})(1-p^{-6\gamma})}=1.
\end{equation}
\begin{equation}\label{E:4.22}
\frac{\sigma_{-\gamma}(6,2;k)} {\sigma_{-\gamma}(5;k)}
= \prod_{p|k} \dfrac{(1-p^{-6\gamma})(1-p^{-2\gamma})}{(1-p^{-\gamma})(1-p^{-5\gamma})}.
\end{equation}

\end{document}